\documentclass[subeqn]{siamltex}

\title{Noncommutative extensions of\\
Ramanujan's ${}_1\psi_1$ summation\thanks{Supported by
Austrian Science Fund FWF,
grant P17563-N13, and partly by EC's IHRP Programme,
grant HPRN-CT-2001-00272 ``Algebraic Combinatorics in Europe''}}

\author{Michael Schlosser\thanks{Fakult\"at f\"ur Mathematik,
Universit\"at Wien, Nordbergstra{\ss}e 15, A-1090 Wien, Austria.\hfill\break
{\it E-mail address:} {\tt schlosse@ap.univie.ac.at}, {\it URL:}
{\tt http://www.mat.univie.ac.at/{\textasciitilde}schlosse}}}

\begin{document}
\maketitle

\begin{abstract}
Using functional equations, we derive noncommutative extensions
of Ramanujan's ${}_1\psi_1$ summation.
\end{abstract}


\section{Introduction}

Hypergeometric series with noncommutative parameters and argument,
in the special case involving square matrices, have been the subject
of recent study, see e.g.\ the papers by Duval and Ovsienko~\cite{DO},
Gr\"unbaum~\cite{G}, Tirao~\cite{T}, and some of the references
mentioned therein. Of course, this subject is also closely related to
the theory of orthogonal matrix polynomials which was initiated by
Krein~\cite{Kr} and has experienced a steady development, see
e.g.\ Dur\'an and L\'opez-Rodr\'{\i}guez~\cite{DL}.

Very recently, Tirao~\cite{T} considered a particular type
of a matrix valued hypergeometric function (which, in our terminology,
belongs to noncommutative hypergeometric series of ``type I'').
He showed, in particular, that the matrix valued hypergeometric
function satisfies the matrix valued hypergeometric differential
equation, and conversely that any solution of the latter is a
matrix valued hypergeometric function. 

In \cite{S2}, the present author investigated hypergeometric and
basic hypergeometric series involving noncommutative parameters and
argument (short: {\em noncommutative hypergeometric series}, and
{\em noncommutative basic} or {\em $Q$-hypergeometric series})
over a unital ring $R$
(or, when considering nonterminating series, over a unital Banach
algebra $R$) from a different, nevertheless completely elementary,
point of view. These investigations were exclusively devoted to the
derivation of summation formulae (which quite surprisingly even exist
in the noncommutative case), aiming to build up a theory of explicit
identities analogous to the rich theory of identities for hypergeometric
and basic hypergeometric series in the classical, commutative case
(cf.\ \cite{Sl} and \cite{GR}). Two closely related types of
noncommmutative series, of ``type I'' and ``type II'', were considered
in \cite{S2}. Most of the summations obtained there concern terminating
series and were proved by induction. An exception are the
noncommutative extensions of the nonterminating $q$-binomial theorem
\cite[Th.~7.2]{S2} which were established using functional equations.
Aside from the latter and some conjectured $Q$-Gau{\ss} summations,
no other explicit summations for nonterminating noncommutative basic
hypergeometric series were given. Furthermore, noncommutative
{\em bilateral} basic hypergeometric series were not even considered.

In this paper, we define noncommutative bilateral basic hypergeometric
series of type I and type II (over an abstract unital Banach algebra $R$)
and prove,  using functional equations, noncommutative extensions of
Ramanujan's ${}_1\psi_1$ summation. These generalize the noncommutative
$Q$-binomial theorem of \cite[Th.~7.2]{S2}. Our proof of the
${}_1\psi_1$ sum here is similar to Andrews and Askey's~\cite{AA} proof
in the classical commutative case. Ramanujan's ${}_1\psi_1$ summation
(displayed in (\ref{11gl})) is one of the fundamental identities in
$q$-series. It is thus just natural to look for different extensions,
including noncommutative ones.

This paper is organized as follows. In Section~\ref{secprel} we review
some standard notations for basic hypergeometric series and then explain
the notation we utilize in the noncommutative case. Section~\ref{sec11}
is devoted to the derivation of noncommutative ${}_1\psi_1$ summations.

We stress again, as in \cite{S2}, that by ``noncommutative'' we do not
mean ``$q$-com\-mutative'' or ``quasi-commutative'' (i.e., where the
variables satisfy a relation like $yx=qxy$; such series are considered
e.g.\ in \cite{K} and \cite{V}) but that the parameters in our series
are elements of some noncommutative unital ring (or unital Banach algebra).

\section{Preliminaries}\label{secprel}

\subsection{Classical (commutative) basic hypergeometric series}

For convenience, we recall some standard notations for basic
hypergeometric series (cf.\ \cite{GR}). When considering the
noncommutative extensions in Subsection~\ref{subsecnc} and in
Section~\ref{sec11}, the reader may find it useful to compare
with the classical, commutative case.

Let $q$ be a complex number such that $0<|q|<1$. Define the
{\em $q$-shifted factorial} for all integers $k$ (including
infinity) by
$$
(a;q)_k:=\prod_{j=1}^k(1-aq^j).
$$
We write
\begin{equation}\label{defhyp}
{}_r\phi_{r-1}\!\left[\!\begin{array}{c}a_1,a_2,\dots,a_r\\
b_1,b_2,\dots,b_{r-1}\end{array}\!;q,z\right]:=
\sum _{k=0}^{\infty}\frac {(a_1;q)_k(a_2;q)_k\dots (a_r;q)_k}
{(q,q)_k(b_1;q)_k\dots(b_{r-1};q)_k}\,z^k,
\end{equation}
to denote the (unilateral)
{\em basic hypergeometric ${}_r\phi_{r-1}$ series}.
Further, we write
\begin{equation}\label{defbhyp}
{}_r\psi_r\!\left[\!\begin{array}{c}a_1,a_2,\dots,a_r\\
b_1,b_2,\dots,b_r\end{array}\!;q,z\right]:=
\sum_{k=-\infty}^{\infty}\frac {(a_1;q)_k\dots (a_r;q)_k}
{(b_1;q)_k\dots(b_r;q)_k}\,z^k,
\end{equation}
to denote the {\em bilateral basic hypergeometric ${}_r\psi_r$ series}.

In (\ref{defhyp}) and (\ref{defbhyp}), $a_1,\dots,a_r$ are
called the {\em upper parameters}, $b_1,\dots,b_r$ the
{\em lower parameters}, $z$ is the {\em argument}, and $q$ the
{\em base} of the series.
The bilateral ${}_r\psi_r$ series in (\ref{defbhyp}) reduces to
a unilateral ${}_r\phi_{r-1}$ series if one of the lower parameters,
say $b_r$, equals $q$ (or more generally, an integral power of $q$).

The basic hypergeometric ${}_r\phi_{r-1}$ series
terminates if one of the upper parameters, say $a_r$, is of the
form $q^{-n}$, for a nonnegative integer $n$.
If the basic hypergeometric series does not terminate then it converges
by the ratio test when $|z|<1$.
Similarly, the bilateral basic hypergeometric series
converges when $|z|<1$ and $|b_1\dots b_r/a_1\dots a_rz|<1$.

We recall two important summations. One of them
is the (nonterminating) $q$-binomial theorem,
\begin{equation}\label{10gl}
{}_1\phi_0\!\left[\!\begin{array}{c}a\\-\end{array}\!;q,z\right]
=\frac{(az;q)_\infty}{(z;q)_\infty},
\end{equation}
where $|z|<1$ (cf.\ \cite[Appendix~(II.3)]{GR}).
It was discovered independently by several
mathematicians, including Cauchy, Gau{\ss}, and Heine.
A bilateral extension of (\ref{10gl}) was found by the legendary
Indian mathematician Ramanujan (see Hardy~\cite{Hr}),
\begin{equation}\label{11gl}
{}_1\psi_1\!\left[\!\begin{array}{c}a\\b\end{array}\!;q,z\right]
=\frac{(q;q)_\infty(b/a;q)_\infty(az;q)_\infty(q/az;q)_\infty}
{(b;q)_\infty(q/a;q)_\infty(z;q)_\infty(b/az;q)_\infty},
\end{equation}
where $|z|<1$ and $|b/az|<1$ (cf.\ \cite[Appendix~(II.29)]{GR}).
Unfortunately, Ramanujan (who very rarely gave any proofs)
did not provide a proof of the above bilateral summation. 
The first proof of (\ref{11gl}) was given by
Hahn~\cite[$\kappa=0$ in Eq.~(4.7)]{H}.
Other proofs were given by Jackson~\cite{J}, Ismail~\cite{I},
Andrews and Askey~\cite{AA}, the author~\cite[Sec.~3]{S1}, and others.
Some immediate applications of Ramanujan's summation formula to
arithmetic number theory are considered in \cite[Sec.~10.6]{AAR}.

\subsection{Noncommutative basic hypergeometric series}\label{subsecnc}

Most of the following definitions are taken from \cite{S2}.
However, the definitions for noncommutative {\em bilateral}
basic hypergeometric series in (\ref{defncbhypIQ}) and
(\ref{defncbhypIIQ}) (although obvious) are new.

Let $R$ be a unital ring (i.e., a ring with a multiplicative identity).
When considering infinite series and infinite products of elements
of $R$ we shall further assume that $R$ is a Banach algebra
(with some norm $\|\cdot\|$). The elements of $R$ will be denoted by
capital letters $A,B,\dots$. In general these elements do not commute
with each other; however, we may sometimes specify certain commutation
relations explicitly. We denote the identity by $I$ and the zero element
by $O$. Whenever a multiplicative inverse element exists for any $A\in R$,
we denote it by $A^{-1}$. (Since $R$ is a unital ring, we have
$AA^{-1}=A^{-1}A=I$.) On the other hand, as we shall implicitly assume
that all the expressions which appear are well defined, whenever
we write $A^{-1}$ we assume its existence.
For instance, in (\ref{defncpochIQ}) and (\ref{defncpochIIQ})
we assume that $I-B_iQ^j$ is invertible
for all $1\le i\le r$, $0\le j<k$.

An important special case is when $R$ is the ring of $n\times n$
square matrices (our notation is certainly suggestive with respect
to this interpretation), or, more generally, one may view
$R$ as a space of some abstract operators.

Let $\mathsf Z$ be the set of integers.
For $l,m\in\mathsf Z\cup\{\pm\infty\}$ we define
the noncommutative product as follows:
\begin{equation}\label{ncprod}
\prod_{j=l}^mA_j=\left\{
\begin{array}{ll}1
&m=l-1\\
A_lA_{l+1}\dots A_m&m\ge l\\
A_{l-1}^{-1}A_{l-2}^{-1}\dots A_{m+1}^{-1}&m<l-1
\end{array}\right..
\end{equation}
Note that
\begin{equation}\label{invprod}
\prod_{j=l}^mA_j=\prod_{j=m+1}^{l-1}A_{m+l-j}^{-1},
\end{equation}
for all $l,m\in\mathsf Z\cup\{\pm\infty\}$.

Throughout this paper, $Q$ will be a parameter which commutes with
any of the other parameters appearing in the series.
(For instance, a central element such as $Q=qI$, a scalar multiple
of the unit element in $R$,
for $qI\in R$, trivially satisfies this requirement.)

Let $k\in\mathsf Z\cup\{\infty\}$.
We define the generalized {\em noncommutative $Q$-shifted factorial
of type I}\/ by
\begin{equation}\label{defncpochIQ}
\left\lceil\!\begin{array}{c}A_1,A_2,\dots,A_r\\
B_1,B_2,\dots,B_r\end{array}\!;Q,Z\right\rfloor_k:=
\prod_{j=1}^k\left[\left(\prod_{i=1}^r
(I-B_iQ^{k-j})^{-1}(I-A_iQ^{k-j})\right)Z\right].
\end{equation}
Similarly, we define the generalized
{\em noncommutative $Q$-shifted factorial of type II}\/ by
\begin{equation}\label{defncpochIIQ}
\left\lfloor\!\begin{array}{c}A_1,A_2,\dots,A_r\\
B_1,B_2,\dots,B_r\end{array}\!;Q,Z\right\rceil_k:=
\prod_{j=1}^k\left[\left(\prod_{i=1}^r
(I-B_iQ^{j-1})^{-1}(I-A_iQ^{j-1})\right)Z\right].
\end{equation}

Note the unusual usage of brackets (``floors'' and ``ceilings''
are intermixed) on the left-hand sides of (\ref{defncpochIQ}) and
(\ref{defncpochIIQ}) which is intended to suggest that the products
involve noncommuting factors in a prescribed order.
In both cases, the product,
read from left to right, starts with a denominator factor.
The brackets in the form ``$\lceil-\rfloor$'' are intended to denote
that the factors are {\em falling},
while in ``$\lfloor-\rceil$'' that they are {\em rising}.

We define the {\em noncommutative basic
hypergeometric series of type I}\/ by
\begin{equation}\label{defnchypIQ}
{}_r\phi_{r-1}\!\left\lceil\!\begin{array}{c}A_1,A_2,\dots,A_r\\
B_1,B_2,\dots,B_{r-1}\end{array}\!;Q,Z\right\rfloor:=
\sum_{k\ge 0}
\left\lceil\!\begin{array}{c}A_1,A_2,\dots,A_r\\
B_1,B_2,\dots,B_{r-1},Q\end{array}\!;Q,Z\right\rfloor_k,
\end{equation}
and the {\em noncommutative basic
hypergeometric series of type II}\/ by
\begin{equation}\label{defnchypIIQ}
{}_r\phi_{r-1}\!\left\lfloor\!\begin{array}{c}A_1,A_2,\dots,A_r\\
B_1,B_2,\dots,B_{r-1}\end{array}\!;Q,Z\right\rceil:=
\sum_{k\ge 0}
\left\lfloor\!\begin{array}{c}A_1,A_2,\dots,A_r\\
B_1,B_2,\dots,B_{r-1},Q\end{array}\!;Q,Z\right\rceil_k.
\end{equation}
Further, we define the {\em noncommutative bilateral basic
hypergeometric series of type I}\/ by
\begin{equation}\label{defncbhypIQ}
{}_r\psi_r\!\left\lceil\!\begin{array}{c}A_1,A_2,\dots,A_r\\
B_1,B_2,\dots,B_r\end{array}\!;Q,Z\right\rfloor:=
\sum_{k=-\infty}^\infty
\left\lceil\!\begin{array}{c}A_1,A_2,\dots,A_r\\
B_1,B_2,\dots,B_r\end{array}\!;Q,Z\right\rfloor_k,
\end{equation}
and the {\em noncommutative bilateral basic
hypergeometric series of type II}\/ by
\begin{equation}\label{defncbhypIIQ}
{}_r\psi_r\!\left\lfloor\!\begin{array}{c}A_1,A_2,\dots,A_r\\
B_1,B_2,\dots,B_r\end{array}\!;Q,Z\right\rceil:=
\sum_{k=-\infty}^\infty
\left\lfloor\!\begin{array}{c}A_1,A_2,\dots,A_r\\
B_1,B_2,\dots,B_r\end{array}\!;Q,Z\right\rceil_k.
\end{equation}

We also refer to the respective series as
{\em (noncommutative) $Q$-hypergeometric series}.
In each case (of type I and type II), the ${}_r\phi_{r-1}$ series
terminates if one of the upper parameters $A_i$ is of the form $Q^{-n}$.
If the ${}_r\phi_{r-1}$ series does not terminate, then
(implicitly assuming that $R$ is a unital Banach algebra with some
norm $\|\cdot\|$) it converges when $\|Z\|<1$.
Similarly, the ${}_r\psi_r$ series converges in $R$ when $\|Z\|<1$ and
$\|Z^{-1}\prod_{i=1}^rA_{r+1-i}^{-1}B_{r+1-i}\|<1$.

We also consider reversed (or ``transposed'') versions of
generalized noncommutative $Q$-shifted factorials and
noncommutative bilateral basic hypergeometric series of type I and II.
These are defined as follows (compare with (\ref{defncpochIQ}),
(\ref{defncpochIIQ}), (\ref{defncbhypIQ}) and (\ref{defncbhypIIQ})):
\begin{equation}\label{defncpochIr}
\phantom{xy}^{\scriptstyle\sim\atop{}}\!\!\left\lceil\!
\begin{array}{c}A_1,A_2,\dots,A_r\\
B_1,B_2,\dots,B_r\end{array}\!;Q,Z\right\rfloor_k:=
\prod_{j=1}^k\left(Z\prod_{i=1}^r(I-A_iQ^{j-1})
(I-B_iQ^{j-1})^{-1}\right),
\end{equation}
\begin{equation}\label{defncpochIIr}
\phantom{xy}^{\scriptstyle\sim\atop{}}\!\!\left\lfloor\!
\begin{array}{c}A_1,A_2,\dots,A_r\\
B_1,B_2,\dots,B_r\end{array}\!;Q,Z\right\rceil_k:=
\prod_{j=1}^k\left(Z\prod_{i=1}^r(I-A_iQ^{k-j})
(I-B_iQ^{k-j})^{-1}\right),
\end{equation}
\begin{equation}\label{defnchypIr}
{}_r\psi_r\!^{\scriptstyle\sim\atop{}}\!\!\left\lceil\!
\begin{array}{c}A_1,A_2,\dots,A_r\\
B_1,B_2,\dots,B_r\end{array}\!;Q,Z\right\rfloor:=
\sum_{k=-\infty}^\infty
{}^{\scriptstyle\sim\atop{}}\!\!\left\lceil\!
\begin{array}{c}A_1,A_2,\dots,A_r\\
B_1,B_2,\dots,B_r\end{array}\!;Q,Z\right\rfloor_k,
\end{equation}
and
\begin{equation}\label{defnchypIIr}
{}_r\psi_r\!^{\scriptstyle\sim\atop{}}\!\!\left\lfloor\!
\begin{array}{c}A_1,A_2,\dots,A_r\\
B_1,B_2,\dots,B_r\end{array}\!;Q,Z\right\rceil:=
\sum_{k=-\infty}^\infty
{}^{\scriptstyle\sim\atop{}}\!\!\left\lfloor\!
\begin{array}{c}A_1,A_2,\dots,A_r\\
B_1,B_2,\dots,B_r\end{array}\!;Q,Z\right\rceil_k.
\end{equation}
Of course, reversed versions of the unilateral noncommutative basic
hypergeometric series are defined analogously.

\section{Noncommutative ${}_1\psi_1$ summations}\label{sec11}

In \cite[Th.~7.2]{S2}, the following two noncommutative extensions of
the nonterminating $q$-binomial theorem (which generalize \cite[II.3]{GR})
were given.

\begin{proposition}\label{nc1f0Q}
Let $A$ and $Z$ be noncommutative parameters of some unital
Banach algebra, and suppose that $Q$ commutes with both $A$ and $Z$.
Further, assume that $\|Z\|<1$. Then we have the following
summation for a noncommutative basic hypergeometric series of type I.
\begin{equation}\label{nc1f0QIgl}
{}_1\phi_0\!\left\lceil\!\begin{array}{c}A\\-\end{array}\!;
Q,Z\right\rfloor=
\left\lfloor\!\begin{array}{c}AZ\\Z\end{array}\!;Q,I\right\rceil_\infty.
\end{equation}
Further, we we have the following
summation for a noncommutative basic hypergeometric series of type II.
\begin{equation}\label{nc1f0QIIgl}
{}_1\phi_0\!\left\lfloor\!\begin{array}{c}A\\-\end{array}\!;Q,Z\right\rceil=
{}^{\scriptstyle\sim\atop{}}\!\!
\left\lfloor\!\begin{array}{c}AZ\\Z\end{array}\!;Q,I\right\rceil_\infty.
\end{equation}
\end{proposition}

Here we extend Proposition~\ref{nc1f0Q} to summations
for bilateral series. Our proof is similar to that of the
classical result given in \cite{AA} (see also
\cite[p.~502, first proof of Th.~10.5.1]{AAR}), but also
similar to the proof of Proposition~\ref{nc1f0Q} given in \cite{S2}.

\begin{theorem}\label{nc11}
Let $A$, $B$ and $Z$ be noncommutative parameters of some unital
Banach algebra, suppose that $Q$ and $B$ both commute with any of the
other parameters. Further, assume that $\|Z\|<1$ and $\|BZ^{-1}A^{-1}\|<1$.
Then we have the following summation for a noncommutative
bilateral basic hypergeometric series of type I.
{\setlength\arraycolsep{2pt}\begin{eqnarray}\label{nc11Igl}
&&{}_1\psi_1\!\left\lceil\!
\begin{array}{c}A\\B\end{array}\!;Q,Z\right\rfloor=\\\nonumber&&
\left\lfloor\!\begin{array}{c}Q,BZ^{-1}A^{-1}Z\\
B,BZ^{-1}A^{-1}\end{array}\!;Q,I\right\rceil_\infty
\left\lceil\!\begin{array}{c}Z^{-1}A^{-1}Q\\
Z^{-1}A^{-1}ZQ\end{array}\!;Q,I\right\rfloor_\infty
\left\lfloor\!\begin{array}{c}AZ\\Z\end{array}\!;Q,I\right\rceil_\infty.
\end{eqnarray}}
Further, we have the following summation for a noncommutative
bilateral basic hypergeometric series of type II.
{\setlength\arraycolsep{2pt}\begin{eqnarray}\label{nc11IIgl}
&&{}_1\psi_1\!\left\lfloor\!
\begin{array}{c}A\\B\end{array}\!;
Q,Z\right\rceil=\\\nonumber
&&{}^{\scriptstyle\sim\atop{}}\!\!
\left\lfloor\!\begin{array}{c}AZ\\Z\end{array}\!;Q,I\right\rceil_\infty
\left\lfloor\!\begin{array}{c}Z^{-1}A^{-1}Q,Q\\
B,A^{-1}Q\end{array}\!;Q,I\right\rceil_\infty
{}^{\scriptstyle\sim\atop{}}\!\!\left\lfloor\!\begin{array}{c}BA^{-1}\\
BZ^{-1}A^{-1}\end{array}\!;Q,I\right\rceil_\infty.
\end{eqnarray}}
\end{theorem}

Clearly, Theorem~\ref{nc11} reduces to Proposition~\ref{nc1f0Q} when $B=Q$.

{\em Proof of Theorem~\ref{nc11}.}
We prove (\ref{nc11Igl}), leaving the proof of (\ref{nc11IIgl}),
which is similar,
to the reader.

Let $f(A,B,Z)$ denote the series on the left-hand side of
(\ref{nc11Igl}). We make use of the two simple identities
\begin{subequations}
{\setlength\arraycolsep{2pt}\begin{eqnarray}\label{fq1}
Z&=&AZQ^k+(I-AQ^k)Z,\\\label{fq2}
I&=&BQ^k+(I-BQ^k),
\end{eqnarray}}
\end{subequations}
to obtain two functional equations for $f$.
We also make use of the simple relation
\begin{equation}\label{fq}
f(A,B,Z)=f(AQ,BQ,Z)(I-B)^{-1}(I-A)Z,
\end{equation}
obtained by shifting the summation index in $f$ by one.

First, (\ref{fq1}) gives
\begin{equation}\label{fq11}
Zf(A,B,Z)=AZf(A,B,ZQ)+f(AQ,B,Z)(I-A)Z,
\end{equation}
while (\ref{fq2}) gives
\begin{equation}\label{fq21}
f(A,BQ,Z)=Bf(A,BQ,ZQ)+(I-B)f(A,B,Z).
\end{equation}
Combining (\ref{fq21}), (\ref{fq11}), and (\ref{fq}), one readily deduces
{\setlength\arraycolsep{2pt}\begin{eqnarray*}
f(A,BQ,Z)&=&(I-B)f(A,B,Z)+BZ^{-1}A^{-1}Zf(A,BQ,Z)\\&&{}
-BZ^{-1}A^{-1}f(AQ,BQ,Z)(I-A)Z\\&=&
(I-B)f(A,B,Z)+BZ^{-1}A^{-1}Zf(A,BQ,Z)\\&&{}
-BZ^{-1}A^{-1}f(A,B,Z)(I-B),
\end{eqnarray*}}
or equivalently
$$
(I-BZ^{-1}A^{-1}Z)f(A,BQ,Z)=(I-B)(I-BZ^{-1}A^{-1})f(A,B,Z),
$$
thus
\begin{equation}\label{fq0}
f(A,B,Z)=(I-BZ^{-1}A^{-1})^{-1}(I-BZ^{-1}A^{-1}Z)(I-B)^{-1}
f(A,BQ,Z).
\end{equation}
Iteration of (\ref{fq0}) gives 
{\setlength\arraycolsep{2pt}\begin{eqnarray}\label{fq01}
&&f(A,B,Z)=\\\nonumber&&\prod_{j=0}^\infty\left[(I-BZ^{-1}A^{-1}Q^j)^{-1}
(I-BZ^{-1}A^{-1}ZQ^j)(I-BQ^j)^{-1}\right]
f(A,O,Z).
\end{eqnarray}}
We still need to compute $f(A,O,Z)$. It is not easy to do this directly
but we know the value of $f(A,Q,Z)$ (by Proposition~\ref{nc1f0Q}).
We set $B=Q$ in (\ref{fq01}) which gives
{\setlength\arraycolsep{2pt}\begin{eqnarray*}
&&f(A,Q,Z)=\\\nonumber&&\prod_{j=0}^\infty\left[(I-Z^{-1}A^{-1}Q^{j+1})^{-1}
(I-Z^{-1}A^{-1}ZQ^{j+1})(I-Q^{j+1})^{-1}\right]
f(A,O,Z),
\end{eqnarray*}}
thus we obtain
\begin{equation}\label{fq03}
f(A,O,Z)=\bigg[\prod_{j=0}^\infty(I-Q^{j+1})\bigg]
\left\lceil\!\begin{array}{c}Z^{-1}A^{-1}Q\\
Z^{-1}A^{-1}ZQ\end{array}\!;Q,I\right\rfloor_\infty
f(A,Q,Z).
\end{equation}
Combination of (\ref{fq01}), (\ref{fq03}) and (\ref{nc1f0QIgl})
establishes the result. \endproof

In the ${}_1\psi_1$ summations of Theorem~\ref{nc11}, the lower parameter
$B$ commutes with both $A$ and $Z$ while $A$ does not commute with $Z$.
In the next theorem the roles of $A$ and $B$ are interchanged. Here $A$
commutes with both $B$ and $Z$ while $B$ does not commute with $Z$.

\begin{theorem}\label{nc11r}
Let $A$, $B$ and $Z$ be noncommutative parameters of some Banach algebra,
suppose that $Q$ and $A$ both commute with any of the other parameters.
Further, assume that $\|Z\|<1$ and $\|BZ^{-1}A^{-1}\|<1$.
Then we have the following summation for a noncommutative
bilateral basic hypergeometric series of type I.
{\setlength\arraycolsep{2pt}\begin{eqnarray}\label{nc11rIgl}
&&{}_1\psi_1\!\left\lceil\!
\begin{array}{c}A\\B\end{array}\!;Q,Z\right\rfloor=\\\nonumber&&
Z^{-1}\left\lfloor\!\begin{array}{c}Q,ZBZ^{-1}A^{-1}\\
A^{-1}Q,Z\end{array}\!;Q,I\right\rceil_\infty
\left\lceil\!\begin{array}{c}AZ\\
ZBZ^{-1}\end{array}\!;Q,I\right\rfloor_\infty
\left\lfloor\!\begin{array}{c}Z^{-1}A^{-1}Q\\
BZ^{-1}A^{-1}\end{array}\!;Q,I\right\rceil_\infty Z.
\end{eqnarray}}
Further, we have the following summation for a noncommutative
bilateral basic hypergeometric series of type II.
{\setlength\arraycolsep{2pt}\begin{eqnarray}\label{nc11rIIgl}
&&{}_1\psi_1\!\left\lfloor\!
\begin{array}{c}A\\B\end{array}\!;
Q,Z\right\rceil=\\\nonumber
&&Z^{-1}{}^{\scriptstyle\sim\atop{}}\!\!
\left\lfloor\!\begin{array}{c}Z^{-1}A^{-1}Q\\
BZ^{-1}A^{-1}\end{array}\!;Q,I\right\rceil_\infty
\left\lfloor\!\begin{array}{c}AZ,Q\\
A^{-1}Q,B\end{array}\!;Q,I\right\rceil_\infty
{}^{\scriptstyle\sim\atop{}}\!\!\left\lfloor\!\begin{array}{c}A^{-1}B\\
Z\end{array}\!;Q,I\right\rceil_\infty Z.
\end{eqnarray}}
\end{theorem}

\begin{proof}
We indicate the derivation of (\ref{nc11rIgl}) from (\ref{nc11Igl}).
(The derivation of (\ref{nc11rIIgl}) from (\ref{nc11IIgl}) is analogous.)
The sum on the left-hand side of (\ref{nc11Igl}) remains unchanged
if the summation index, say $k$, is replaced by $-k$. Using (\ref{invprod}),
we compute
{\setlength\arraycolsep{2pt}\begin{eqnarray}
\left\lceil\!
\begin{array}{c}A\\B\end{array}\!;Q,Z\right\rfloor_{-k}&=&
\prod_{j=1}^{-k}\left[
(I-BQ^{-k-j})^{-1}(I-AQ^{-k-j})Z\right]\\\nonumber
&=&\prod_{j=1-k}^0\left[(I-BQ^{-1+j})^{-1}
(I-AQ^{-1+j})Z\right]^{-1}\\\nonumber
&=&\prod_{j=1}^k\left[Z^{-1}(I-AQ^{-1-k+j})^{-1}
(I-BQ^{-1-k+j})\right]\\\nonumber
&=&\prod_{j=1}^k\left[Z^{-1}A^{-1}(I-A^{-1}Q^{1+k-j})^{-1}
(I-B^{-1}Q^{1+k-j})B\right]\\\nonumber
&=&Z^{-1}A^{-1}\left\lceil\!
\begin{array}{c}B^{-1}Q\\
A^{-1}Q\end{array}\!;Q,BZ^{-1}A^{-1}\right\rfloor_k AZ.
\end{eqnarray}}
Thus, by performing the simultaneous replacements $A\mapsto B^{-1}Q$,
$B\mapsto A^{-1}Q$, $Z\mapsto BZ^{-1}A^{-1}$, in (\ref{nc11Igl}),
we obtain (\ref{nc11rIgl}).
\end{proof}

We complete this paper with four more ${}_1\psi_1$ summations,
immediately obtained from corresponding summations in
Theorems~\ref{nc11} and \ref{nc11r} by ``reversing all products''
(cf.\ \cite[Subsec.~8.2]{S2}) on each side of the respective identities.

\begin{theorem}\label{nc11t}
Let $A$, $B$ and $Z$ be noncommutative parameters of some unital
Banach algebra, suppose that $Q$ and $B$ both commute with any of the
other parameters. Further, assume that $\|Z\|<1$ and $\|A^{-1}Z^{-1}B\|<1$.
Then we have the following summation for a reversed noncommutative
bilateral basic hypergeometric series of type I.
{\setlength\arraycolsep{2pt}\begin{eqnarray}\label{nc11tIgl}
&&{}_1\psi_1^{\scriptstyle\sim\atop{}}\!\!\left\lceil\!
\begin{array}{c}A\\B\end{array}\!;Q,Z\right\rfloor=\\\nonumber&&
{}^{\scriptstyle\sim\atop{}}\!\!
\left\lfloor\!\begin{array}{c}ZA\\Z\end{array}\!;Q,I\right\rceil_\infty
{}^{\scriptstyle\sim\atop{}}\!\!\left\lceil\!\begin{array}{c}A^{-1}Z^{-1}Q\\
ZA^{-1}Z^{-1}Q\end{array}\!;Q,I\right\rfloor_\infty
{}^{\scriptstyle\sim\atop{}}\!\!
\left\lfloor\!\begin{array}{c}ZA^{-1}Z^{-1}B,Q\\
A^{-1}Z^{-1}B,B\end{array}\!;Q,I\right\rceil_\infty.
\end{eqnarray}}
Further, we have the following summation for a noncommutative
bilateral basic hypergeometric series of type II.
{\setlength\arraycolsep{2pt}\begin{eqnarray}\label{nc11tIIgl}
&&{}_1\psi_1^{\scriptstyle\sim\atop{}}\!\!\left\lfloor\!
\begin{array}{c}A\\B\end{array}\!;
Q,Z\right\rceil=\\\nonumber
&&\left\lfloor\!\begin{array}{c}A^{-1}B\\
A^{-1}Z^{-1}B\end{array}\!;Q,I\right\rceil_\infty
{}^{\scriptstyle\sim\atop{}}\!\!
\left\lfloor\!\begin{array}{c}Q,A^{-1}Z^{-1}Q\\
A^{-1}Q,B\end{array}\!;Q,I\right\rceil_\infty
\left\lfloor\!\begin{array}{c}ZA\\Z\end{array}\!;Q,I\right\rceil_\infty.
\end{eqnarray}}
\end{theorem}

\begin{theorem}\label{nc11tr}
Let $A$, $B$ and $Z$ be noncommutative parameters of some Banach algebra,
suppose that $Q$ and $A$ both commute with any of the other parameters.
Further, assume that $\|Z\|<1$ and $\|A^{-1}Z^{-1}B\|<1$.
Then we have the following summation for a noncommutative
bilateral basic hypergeometric series of type I.
{\setlength\arraycolsep{2pt}\begin{eqnarray}\label{nc11trIgl}
\phantom{xyz}&&{}_1\psi_1^{\scriptstyle\sim\atop{}}\!\!\left\lceil\!
\begin{array}{c}A\\B\end{array}\!;Q,Z\right\rfloor=\\\nonumber&&
Z{}^{\scriptstyle\sim\atop{}}\!\!
\left\lfloor\!\begin{array}{c}A^{-1}Z^{-1}Q\\
A^{-1}Z^{-1}B\end{array}\!;Q,I\right\rceil_\infty
{}^{\scriptstyle\sim\atop{}}\!\!\left\lceil\!\begin{array}{c}ZA\\
Z^{-1}BZ\end{array}\!;Q,I\right\rfloor_\infty
{}^{\scriptstyle\sim\atop{}}\!\!
\left\lfloor\!\begin{array}{c}A^{-1}Z^{-1}BZ,Q\\
Z,A^{-1}Q\end{array}\!;Q,I\right\rceil_\infty Z^{-1}.
\end{eqnarray}}
Further, we have the following summation for a noncommutative
bilateral basic hypergeometric series of type II.
{\setlength\arraycolsep{2pt}\begin{eqnarray}\label{nc11trIIgl}
&&{}_1\psi_1^{\scriptstyle\sim\atop{}}\!\!\left\lfloor\!
\begin{array}{c}A\\B\end{array}\!;
Q,Z\right\rceil=\\\nonumber
&&Z\left\lfloor\!\begin{array}{c}BA^{-1}\\
Z\end{array}\!;Q,I\right\rceil_\infty
{}^{\scriptstyle\sim\atop{}}\!\!\left\lfloor\!\begin{array}{c}Q,ZA\\
B,A^{-1}Q\end{array}\!;Q,I\right\rceil_\infty
\left\lfloor\!\begin{array}{c}A^{-1}Z^{-1}Q\\
A^{-1}Z^{-1}B\end{array}\!;Q,I\right\rceil_\infty
Z^{-1}.
\end{eqnarray}}
\end{theorem}

\end{document}